\documentclass[11pt]{article}

\usepackage{IEEEtrantools}
\usepackage{latexsym}
\usepackage{amsfonts}
\usepackage{amssymb}
\usepackage{amsmath,amsfonts,amsthm}
\usepackage{mathrsfs}
\usepackage{multirow}
\usepackage{tikz}
 \usetikzlibrary{arrows,decorations.markings}

\newcommand{\be}{\begin{equation}}
\newcommand{\ee}{\end{equation}}
\newcommand{\bea}{\begin{eqnarray}}
\newcommand{\eea}{\end{eqnarray}}
\newcommand{\bean}{\begin{eqnarray*}}
\newcommand{\eean}{\end{eqnarray*}}
\newcommand{\brray}{\begin{array}}
\newcommand{\erray}{\end{array}}
\newcommand{\biearray}{\begin{IEEEarray}{rCl}}
\newcommand{\eiearray}{\end{IEEEarray}}

\newcommand{\newsection}[1]{\setcounter{equation}{0}
\setcounter{dfn}{0}
\section{#1}}

\newtheorem{dfn}{Definition}[section]
\newtheorem{thm}[dfn]{Theorem}
\newtheorem{lmma}[dfn]{Lemma}
\newtheorem{ppsn}[dfn]{Proposition}
\newtheorem{crlre}[dfn]{Corollary}
\newtheorem{xmpl}[dfn]{Example}
\newtheorem{rmrk}[dfn]{Remark}

\newcommand{\bdfn}{\begin{dfn}\rm}
\newcommand{\bthm}{\begin{thm}}
\newcommand{\blmma}{\begin{lmma}}
\newcommand{\bppsn}{\begin{ppsn}}
\newcommand{\bcrlre}{\begin{crlre}}
\newcommand{\bxmpl}{\begin{xmpl}}
\newcommand{\brmrk}{\begin{rmrk}\rm}

\newcommand{\edfn}{\end{dfn}}
\newcommand{\ethm}{\end{thm}}
\newcommand{\elmma}{\end{lmma}}
\newcommand{\eppsn}{\end{ppsn}}
\newcommand{\ecrlre}{\end{crlre}}
\newcommand{\exmpl}{\end{xmpl}}
\newcommand{\ermrk}{\end{rmrk}}

\newcommand{\bbc}{\mathbb{C}}
\newcommand{\bbz}{\mathbb{Z}}
\newcommand{\bbn}{\mathbb{N}}
\newcommand{\bbr}{\mathbb{R}}

\newcommand{\cla}{\mathcal{A}}

\newcommand{\clh}{\mathcal{H}}

\newcommand{\cll}{\mathcal{L}}

\newcommand{\cls}{\mathcal{S}}

\newcommand{\prf}{\noindent{\it Proof\/}: }

\def \qed { \mbox{}\hfill
$\Box$\vspace{1ex}}

\begin{document}

 \author{\sc  { Bipul Saurabh}}
\title{ Spectral dimension of  spheres }
\maketitle
 \begin{abstract}
In this paper, we associate a growth graph and a  length operator  to  a quotient space of a semisimple compact Lie group. Under certain assumptions, we 
show that the spectral dimension of a homogeneous space is greater than or equal to  summability of the length operator. Using this,
we compute spectral dimensions of spheres. 
\end{abstract}
{\bf AMS Subject Classification No.:} {\large 58}B{\large 34}, {\large 58}B{\large 32}, {\large46}L{\large 87} 

\newsection{Introduction}  

Motivated by Connes' definition of  dimension of a spectral triple,   Chakraborty and Pal
(\cite{ChaPal-2016aa}) introduced an invariant called spectral dimension, for an ergodic $C^*$-dynamical system or equivalently 
for a homogeneous space of a compact quantum group.  Ergodicity of a $C^*$-dynamical system gives a unique 
state invariant under group action.  Chakraborty and Pal considered 
all finitely summable equivariant  spectral triples on the GNS space of the invariant  state  
and defined the spectral dimension of the homogeneous space  to be the infimum of
the summability of the associated Dirac operators. Some questions naturally arise about this invariant.
\begin{enumerate}
 \item  For  classical homogeneous spaces, is it some  known
quantity associated with the space?
\item
Is spectral dimension of  $q$-deformation $G_q$ of a semisimple simpy connected compact Lie group $G$ 
is same as $G$?
\item 
Given a Poisson Lie subgroup $H$ of $G$, is spectral dimension of $G_q/H_q$ equal to that of $G/H$? 
\end{enumerate}
To answer all these questions, we need to compute the invariant in several cases. 
Chakraborty and Pal computed 
spectral dimension of many  homogeneous spaces, both in classical and quantum situations and  it was conjectured 
that the spectral dimension of a homogeneous space of a (classical) compact Lie group 
is same as its dimension as a differentiable manifold. 
The spectral dimension of  $SU(2)$ points
towards this conjecture.  Saurabh \cite{Sau-2018aa} considered the quaternion sphere $SP(2n)/SP(2n-2)$ with canonical $SP(2n)$ action and proved that 
its spectral dimension is equal to the dimension of quaternion spheres as a real manifold which strengthens the conjecture of Chakraborty and Pal. 
So far, we have only these two instances which one can cite as a strong evidence in support of the conjecture.
To  gain more confidence about proving the conjecture, it is important to explore more examples. 
In this article, we  take up  
three cases.  First one is odd dimensional sphere $SU(n+1)/SU(n)$ with natural $SU(n+1)$ action. It is a type $A$ homogeneous space.  Second one 
is  even dimensional sphere $SO(2n+1)/SO(2n)$ with  canonical $SO(2n+1)$ action. It is an example of type $B$ homogeneous spaces. And third one is 
the type $D$ homogeneous space $SO(2n)/SO(2n-1)$ with canonical $SO(2n)$ action. 
At this moment, it must be pointed out that though $SU(n+1)/SU(n)$ and $SO(2n+2)/SO(2n+1)$ are homeomorphic as a topological space,  
one can not conclude that their spectral dimensions are equal because the $C^*$-dynamical system in both cases 
are not  isomorphic. Moreover, for $q \in (0,1)$, the spectral dimension of $SU_q(n+1)/SU_q(n)$ is computed in \cite{ChaPal-2016aa} but for $q=1$,
one can not follow that method as 
the behaviour of these classical spaces can be quite different from their  quantum analog. Therefore
it is worthwhile to compute the invariant in these cases. 
In this paper,  we prove that  the spectral dimensions in all the three cases are  equal to their dimensions as a real manifold. So, our results support 
the conjecture of Chakraborty and Pal.

For a subset $S$ of a $C^*$-algebra $A$,   $\overline{S}$ is the closed linear span of $S$ in $A$. We will sometimes write a spectral triple $(\cla,\clh,D)$ 
as $(\clh,\pi,D)$ where $\pi$ is the representation of $\cla$ in the Hilbert space $\clh$.   We denote by $\cls^n$ the $n$-dimensional sphere.

\newsection{Spectral dimension}
In this section, we recall from \cite{ChaPal-2016aa} the definition of spectral dimension of a $C^*$-dynamical system and then give some conditions on a
homogeneous space to get a lower bound on its spectral dimension. Let us begin with the definition of a homogeneous space. 
\bdfn A compact quantum group $G$ acts on a $C^*$-algebra $A$ if there exists a $*$-homomorphism $\tau: A \rightarrow A \otimes C(G)$ such that 
\begin{enumerate}
 \item $(\tau \otimes id)\tau=(id\otimes \Delta)\tau$,
 \item$\overline{\{(I\otimes b)\tau(a): a\in A, b \in C(G)\}}=A\otimes C(G)$.
\end{enumerate}
where $\Delta$ is the comultiplication map of $G$.
An action $\tau$ is called homogeneous or ergodic if the fixed point subalgebra $\{a \in A : \tau(a)=a\otimes I\}$ is $\bbc I$. 
In that case, the triple $(A,G,\tau)$ is called an ergodic $C^*$-dynamical system and the associated $C^*$-algebra $A$ is called a homogeneous space 
of $G$.
\edfn
A covariant representation of a $C^*$-dynamical system $(A,G,\tau)$ is a pair
$(\pi,U)$ consisting of a representation $\pi:A \rightarrow \cll(\clh)$ and a unitary representation 
of $G$ on $\clh$ such that for all $a \in A$, one has 
\[
 (\pi \otimes id)\tau(a)= U(\pi(a) \otimes I)U^*.
\]
\bdfn
Let $(\pi,U)$ be a covariant representation of a $C^*$-dynamical system $(A,G,\tau)$ and $(\clh,\pi,D)$ be a spectral triple 
for a dense $*$-subalgebra $\cla$ of $A$. We say the spectral triple  $(\clh,\pi,D)$ equivariant with respect
to $(\pi,U)$ if $D\otimes I$ commutes with $U$.
\edfn

Given a homogeneous action $\tau$ of $G$, there is a unique invariant state $\rho$ on the homogeneous space $A$ that satisfies 
\[
 (\rho \otimes id)\tau(a)= \rho(a)I, \qquad a\in A.
\]
Consider the GNS representation $(\clh_{\rho}, \pi_{\rho}, \eta_{\rho})$ of $A$ associated with the state $\rho$. 
Using the invariance property of $\tau$, one can show that the 
action $\tau$ induces a unitary representation $U_{\tau}$ of $G$ on $\clh_{\rho}$ and the pair $(\pi_{\rho}, U_{\tau})$ is a 
covariant representation of the
system $(A,G,\tau)$. Let $\mathcal{O}(G)$ be the dense $*$-Hopf subalgebra of $C(G)$ generated by matrix entries of 
irreducible unitary representations 
of $G$. Define 
\[
 \cla:= \{ a \in A: \tau(a) \in A \otimes_{alg}\mathcal{O}(G)\}.
\]
 From part $(1)$ of Theorem $1.5$ in \cite{Pod-1995aa}, it follows  that $\cla$ is a dense $*$-subalgebra of $A$. Let
 $\xi$  be the class of spectral triples of $\cla$ equivariant
 with respect to the covariant representation $(\pi_{\rho}, U_{\tau})$. The spectral dimension denoted by $\cls dim(A,G,\tau)$ of 
 the $C^*$-dynamical system $(A,G,\tau)$ is defined as follows.
 \[
  \cls dim(A,G,\tau):= \inf\{p>0: \exists  D \mbox{ such that }  (\cla,\clh_{\rho},D) \in \xi \mbox{ and } D \mbox{ is $p$-summable}\}.
 \]
Let $G$ be a semisimple compact Lie group and $H$ be a closed Poisson Lie subgroup of $G$. Let $\phi : C(G) \rightarrow C(H)$ be 
a $C^{*}$-epimorphism  obeying $\Delta\phi=(\phi\otimes \phi)\Delta$ where $\Delta$ is the co-multiplication map of $C(G))$.
In such a case, one defines the quotient space $C(G/H)$ by,
\[
      C(G/H)) = \left\{a\in C(G) : (\phi\otimes id)\Delta(a) = I\otimes a\right\}. 
\] 
Consider the following $G$-action on the quotient space $G/H$:
\begin{IEEEeqnarray}{rCl}
 \tau: C(G/H) &\longrightarrow & C(G/H)\otimes C(G) \nonumber \\
 a &\longmapsto & \Delta a \nonumber.
\end{IEEEeqnarray}
where $\Delta$ is the co-multiplication map of the compact  group  $G$. By theorem $1.5$ of \cite{Pod-1995aa}, we get,
\begin{IEEEeqnarray}{rCl} \label{decom}
C(G/H)= \overline{\oplus_{\lambda \in \widehat{G}}\oplus_{i \in I_{\lambda} }W_{(\lambda,i)}}
\end{IEEEeqnarray}
where $\lambda$ represents the highest weight of  a finite-dimensional irreducible co-representation $u_{\lambda}$ of $C(G)$,
$I_{\lambda}$ is  the multiplicity of $u_{\lambda}$ 
and $ W_{(\lambda, i)}$ corresponds to 
$u_{\lambda}$ in the sense of Podles (see page $4$, \cite{Pod-1995aa}) for all $i \in I_{\lambda}$. 
We will reparametrize the index appearing in the equation  (\ref{decom}) as follows.
\[
 \Gamma:=\{(\lambda, j): \lambda \in \widehat{G}, 1 \leq j \leq I_{\lambda}\}.
\]
Therefore
\begin{IEEEeqnarray}{rCl}
C(G/H)= \overline{\oplus_{\gamma \in \Gamma}W_{\gamma}}. \nonumber
\end{IEEEeqnarray}
We will denote by $N_{\gamma}$ the dimension of $W_{\gamma}$.   Define 
\[
 \mathcal{O}(G/H):= \oplus_{\gamma \in \Gamma}W_{\gamma}. 
\]
Then $\mathcal{O}(G/H)$ is a dense  Hopf $*$-algebra consisting of all $a \in  C(G/H)$ such that $\tau(a)\in  C(G/H)\otimes_{alg}\mathcal{O}(G)$. 
It is not difficult to verify that the system $(C(G/H), G, \tau)$ is an ergodic $C^*$-dynamical system. 
Assume that $\rho$ is the invariant state of $\tau$ and $(\clh_{\rho}, \pi_{\rho}, \eta_{\rho})$ be the associated GNS representation  of $C(G/H)$.
The Hilbert space  $\clh_{\rho}$ has a  basis of the form 
\[
 \{ e_{(\gamma,i)}:=\eta_{\rho}(a_{(\gamma,i)}): \{a_{(\gamma,i)}: 1 \leq i \leq N_{\gamma}\} \mbox{ is a basis of } W_{\gamma}, \gamma \in \Gamma\}.
\]
Let $R=\{r_1,r_2,\cdots, r_k\}\subset \mathcal{O}(G/H)$ and $c>0$.  Define a directed graph $\mathcal{G}_R^c$ as follows.
Take the vertex set  to be $\Gamma$. We write  $\gamma \rightsquigarrow_{r_j} \gamma^{'}$ if  $e_{(\gamma^{'},i^{'})}=r_je_{(\gamma,i)}$ and $\frac{\|  e_{(\gamma,i)}\|}{\|e_{(\gamma^{'},i^{'})}\|}<c$
for some $ 1\leq i \leq N_{\gamma}$ and $ 1\leq i^{'} \leq N_{\gamma^{'}}$.
Define edge set of $\mathcal{G}_R^c$ to be 
\[
 E:=\big\{(\gamma,\gamma^{'}): \gamma \rightsquigarrow_{r_i} \gamma^{'}   \mbox{ for some } 1 \leq i \leq k \big\}.
\]
We write $\gamma \rightarrow \gamma^{'}$ if $(\gamma,\gamma^{'})\in E$. We call the directed graph $\mathcal{G}_R^c=(\Gamma , E)$ a growth graph of $\mathcal{O}(G/H)$. 
We say that the graph $\mathcal{G}_R^c$ has a root if there exists a vertex $\gamma_0$
such that for any $\gamma \in \Gamma$, there is a directed path from $\gamma_0$ to $\gamma$. The vertex $\gamma_0$ will
be called a root of the graph $\mathcal{G}_R^c$.
In such a case, define a length function $\ell_{\gamma_0}:\Gamma \longrightarrow \bbn $ as follows;
\begin{IEEEeqnarray*}{rCl}
 \ell_{\gamma_0}(\gamma)= \begin{cases}
                           1 & \mbox{ if } \gamma =\gamma_0,\\
                           \mbox{ length of a shortest path from } \gamma_0 \mbox{ to } \gamma, & \mbox{ otherwise.} \\
                          \end{cases}
\end{IEEEeqnarray*}
Let $L_{\gamma_0}$ be the unbounded positive operator on $\clh_{\rho}$ with dense domain $ \mathcal{O}(G/H)$ sending  $e_{(\gamma,i)}$ to $\ell_{\gamma_0}(\gamma)e_{(\gamma,i)}$. We call 
$L_{\gamma_0}$ the length operator associated with the length function $\ell_{\gamma_0}$ or the root $\gamma_0$.
The following proposition says that if the graph 
$\mathcal{G}_R^c$ has a root then   the growth of eigenvalues of   Dirac operators which have bounded commutators 
with the elements of the algebra $\mathcal{O}(G/H)$  is less than or equal to the growth of eigenvalues of the  length operator 
associated with the root.
\bppsn \label{bounded growth3}
Let $D: e_{\gamma,i} \mapsto d_{\gamma}e_{\gamma,i}$ be a selfadjoint unbounded operator with compact resolvent acting on the Hilbert space  $\clh_{\rho}$ such that the triple
$(\mathcal{O}(G/H), \clh_{\rho}, D)$ is a spectral triple.  Moreover, assume that there exists a finite set $R$ in  $\mathcal{O}(G/H)$ and $c>0$ such that 
the  graph  $\mathcal{G}_R^c$ has a root $\gamma_0$. Then we have 
\[
 |d_{\gamma}|=O(\ell_{\gamma_{0}}(\gamma)).
\]
\eppsn
\prf Since the Dirac operator $D$ has a bounded commutator with elements of $\mathcal{O}(G/H)$,  we can define $M:= \max\{\|[D,r_j]\|: 1 \leq j \leq k\}$.
Let $(\gamma, \gamma^{'}) \in E$. Then 
there exist  $ 1\leq i \leq N_{\gamma}$, $ 1\leq i^{'} \leq N_{\gamma^{'}}$ and $ 1\leq j \leq k$ such that  $e_{(\gamma^{'},i^{'})}=r_je_{(\gamma,i)}$ and $\frac{\|  e_{(\gamma,i)}\|}{\|e_{(\gamma^{'},i^{'})}\|}<c$. 
Hence we have
\begin{IEEEeqnarray}{rCl}
 \|[D,r_j]e_{(\gamma,i)}\|=\|Dr_je_{(\gamma,i)}-r_jDe_{(\gamma,i)}\|=\|De_{(\gamma^{'},i^{'})}-d_{\gamma}r_je_{(\gamma,i)}\|=|d_{\gamma{'}}-d_{\gamma}|\|e_{(\gamma^{'},i^{'})}\| \nonumber 
\end{IEEEeqnarray}
Therefore
\[
 |d_{\gamma{'}}-d_{\gamma}| = \frac{\|[D,r_j]e_{(\gamma,i)}\|}{\|e_{(\gamma^{'},i^{'})}|}\leq \frac{\|[D,r_j]\|\|e_{(\gamma,i)}\|}{\|e_{(\gamma^{'},i^{'})}|}\leq cM
\]
Take $\gamma \in \Gamma$. Let $\gamma_{0}\rightarrow \gamma_1 \rightarrow \cdots  \rightarrow 
\gamma_{\ell_{\gamma_{0}}(\gamma)}=\gamma$  be a shortest path from $\gamma_0$ to $\gamma$.
\[
 |d_{\gamma}| \leq |d_{\gamma_0}|+ |d_{\gamma_1}-d_{\gamma_0}|+ |d_{\gamma_2}-d_{\gamma_1}| \cdots +
 |d_{\gamma}-d_{\gamma_{(\ell_{\gamma_{0}}(\gamma)-1)}}| \leq cM \ell_{\gamma_{0}}(\gamma)
\]
This proves the claim.
\qed

\bppsn \label{geq}
Assume that there exists a finite set $R$ in  $\mathcal{O}(G/H)$ and $c>0$ such that 
the  graph  $\mathcal{G}_R^c$ has a root $\gamma_0$.   Define
$l=\inf\{p: \mbox{ Tr}(L_{\gamma_0}^{-p})< \infty\}$.                                                                                                                                                                               
Then one has
\[
\cls dim(C(G/H),G,\tau) \geq l.
\]
\eppsn
\prf
Let $(\mathcal{O}(G/H),\clh_\rho,D)$ be an equivariant spectral triple of the system $(C(G/H), G, \tau)$. Then following the
arguments in propositions 5.1-5.3 leading to the statement (5.22)  in  \cite{ChaPal-2016aa}, we can assume that $D$ must be of the form 
\[
 De_{(\gamma,i)}=d_{\gamma}e_{(\gamma,i)}, \qquad i  \in \{1,2,\cdots N_{\gamma}\}, \gamma \in \Gamma.  
\]
By Proposition \ref{bounded growth3}, we get  $|d_{\gamma}|=O(\ell_{\gamma_{0}}(\gamma))$. Therefore $D$ is $p$-summable for any $p>l$. This completes the proof.
\qed


\newsection{$SU(n+1)$ action on $\cls^{2n+1}$}
In this section, we will take $G$ to be  $SU(n+1)$ and $H$ to be $  SU(n)$.
For $1 \leq i,j \leq n+1$, define a continuous map 
\begin{IEEEeqnarray}{rCl}
 u_j^i:SU(n+1) & \rightarrow & \bbc;  \quad A \mapsto a_j^i \nonumber
\end{IEEEeqnarray} 
where $a_j^i$ is the $ij^{th}$ entry of $A \in SU(n+1)$. 
The $C^*$-algebra  $C(SU(n+1))$ is generated by elements of the set $\{u_{j}^{i}:1\leq i,j\leq n+1\}$. 
In the same way,  define the generators  $\{v_{j}^{i}:1\leq i,j\leq n\}$  of $C(SU(n))$. 
Define the map $\Phi: C(SU(n+1)) \rightarrow C(SU(n))$ 
as follows. 
\[ 
\Phi(u_{j}^{i})=\begin{cases} 
                             v_{j}^{i}, 
                             & \mbox{ if } i \neq n+1 \mbox{ or } j \neq n+1, \cr
			 \delta_{ij}, & \mbox{ otherwise. } \cr
			\end{cases}
\] 
The quotient space $SU(n+1)/SU(n)$ can be realized as the $2n+1$-dimensional sphere $\cls^{2n+1}$.  Also, each of 
the generators $\{u_{j}^{n+1}:   1\leq j\leq n+1\}$ can be viewed as 
projection on to a fixed complex coordinate of a point in $\cls^{2n+1}\subset \bbc^{n+1}$.  To describe the set of highest weights of all finite-dimensional irreducible co-representation of $C(SU(n+1))$, define 
\[
 X=\{(\lambda_1, \lambda_2, \cdots ,\lambda_{n+1}):\lambda_1\geq \lambda_2 \geq \cdots \geq \lambda_{n+1}, \lambda_i \in \bbz \mbox { for all } 1\leq i \leq n+1\}
\]
we say that two tuple $(\lambda_1, \lambda_2, \cdots ,\lambda_{n+1})$ and  $(\lambda_1^{'}, \lambda_2^{'}, \cdots ,\lambda_{n+1}^{'})$ are equivalent $(\sim)$ if
$\lambda_1-\lambda_1^{'}=\lambda_2-\lambda_2^{'}= \cdots =\lambda_{n+1}-\lambda_{n+1}^{'}$. Define $\varLambda:=X/\sim$. Then $\varLambda$ 
is the  set of highest weights of all finite-dimensional irreducible co-representation of $C(SU(n+1))$.
Using Zhelobenko branching rule (see Theorem 9, page 74, \cite{Zhe-1962aa} 
and Theorem 1.7 in \cite{Pod-1995aa}), we get  
\begin{displaymath}
I_{\lambda}= \begin{cases}
                            1 & \quad \mbox{if} \quad \lambda_1\geq 0, \lambda_i=0 \quad \mbox{for all}\quad 2 \leq  i \leq n  \mbox{ and } \lambda_{n+1}\leq 0, \cr
                             0 & \quad \mbox{otherwise}. \cr
                            \end{cases}
\end{displaymath}  
We will now find a highest weight vector  for each irreducible co-representation of highest weight  $(\lambda_1,0, \cdots,0, \lambda_{n+1})$ 
which belongs to $C(SU(n+1)/SU(n))$. 
Let  $U(\mathfrak{su}(n))$ be the universal enveloping algebra of the  Lie algebra $\mathfrak{su}(n)$. 
We will view $\mathfrak{su}(n)$ as a subset of $U(\mathfrak{su}(n))$. Then $U(\mathfrak{su}(n))$ is 
generated by  $H_i,E_i, F_i \in  \mathfrak{su}(n)$,  $i=1,2,\cdots,n$, 
satisfying the relations given in  page 160, \cite{KliSch-1997aa}. 
Hopf *-structure of  $U(\mathfrak{su}(n))$ comes from the following maps (see page $18$ and page $21$ of \cite{KliSch-1997aa}):
\[
\Delta(r)=r\otimes 1 +1 \otimes r, \quad S(r)=-r, \quad \epsilon(r)=0, \quad r=r^* \quad  \forall r \in \mathfrak{su}(n).
\]
Denote by $T_1$ the finite dimensional irreducible representation of $U(\mathfrak{su}(n))$ 
with highest weight $(1,0,\cdots ,0)$.
There exists unique nondegenerate dual pairing $\left\langle\cdot,\cdot\right\rangle$ between the
Hopf $*$-algebras $U(\mathfrak{su}(n))$ and $\mathcal{O}(SU(n+1)/SU(n))$ such that
\[
  \left\langle f,u_{l}^{k}\right\rangle = t_{kl}(f);        \hspace{1in} \mbox{for }  k= n+1 \mbox{ and }  1\leq l\leq n+1,\\
\]
where $t_{kl}$ is the matrix element of $T_{1}$. Using this, one can give  the algebra $\mathcal{O}(SU(n+1)/SU(n))$ a $U(\mathfrak{su}(n))$-module
structure in the following way.
\[
 f(a)=(1\otimes \langle f,.\rangle)\Delta a,
\]
where $f \in U(\mathfrak{su}(n))$ and $a \in \mathcal{O}(SU(n+1)/SU(n))$.  
We call an element $b \in \mathcal{O}(SU(n+1)/SU(n))$ a highest weight vector  with highest weight $(\lambda_1,0, \cdots,0, \lambda_{n+1})$ if 
\[
 H_1(b)=\lambda_1b, \quad  \quad H_i(b)= 0  \mbox{ for} \quad 2\leq i\leq n-1, \quad H_n(b)=-\lambda_{n+1}b
\]
and
\[
 E_i(b)=0 \qquad  \mbox{ for } 1\leq i \leq n.
\]
The following proposition describes a   highest weight vector with highest weight $(\lambda_1,0, \cdots,0, \lambda_{n+1})$ explicitly.
\bppsn \label{new}
Let $y= u_{1}^{n+1}$ and $z=(u_{n+1}^{n+1})^*$. For $\lambda_1, \lambda_{n+1} \in \bbn$, define $b^{(\lambda_1,\lambda_{n+1})}=y^{\lambda_1}z^{\lambda_{n+1}}$. 
Then $b^{(\lambda_1,\lambda_{n+1})}$ is a  highest weight vectors  in the algebra
$\mathcal{O}(SU(n+1)/SU(n))$ with highest weight $(\lambda_1,0, \cdots,0, \lambda_{n+1})$. 
\eppsn
\prf  It is not difficult to see that
\begin{displaymath}
E_i(y)=E_i(z)=0 \quad \mbox{for}\quad  1\leq i\leq n.
\end{displaymath}
 Further
 \begin{displaymath}
 H_1(y)=y \mbox{ and } H_i(y)=0, \mbox{ for } i >1,
 \end{displaymath}
 and
 \begin{displaymath}
  H_n(z)=-z \mbox{ and } H_i(z)=0, \mbox{ for } i <n.
 \end{displaymath}
Using this and properties of 
Hopf $*$ algebra pairing (see page $21$ of \cite{KliSch-1997aa}), one can check that 
$b^{(\lambda_1,\lambda_{n+1})}$ is a  highest weight vectors 
with highest weight $(\lambda_1,0, \cdots,0, \lambda_{n+1})$. 
\qed \\
For convenience, we will reparametrize the index. Define 
\[
 \Gamma=\{ (\gamma_1,\gamma_2): \gamma_1,\gamma_2 \in \bbn\}.
\]
Hence we have 
\[
 \mathcal{O}(G/H):= \oplus_{\gamma \in \Gamma}W_{\gamma}. 
\]
where $W_{\gamma}$  corresponds to 
$u_{(\gamma_1,0, \cdots , 0, -\gamma_2)}$ in the sense of Podles (see page $4$, \cite{Pod-1995aa}). 
The invariant state $\rho$ of $\tau$ is the  faithful Haar state $h$ of 
$C(SU(n+1))$ restricted to $C(SU(n+1)/SU(n))$. 
Hence  $\clh_{\rho}$ has a  basis of the form 
\[
 \{ e_{(\gamma,i)}: \{e_{(\gamma,i)}: 1 \leq i \leq N_{\gamma}\} \mbox{ is a basis of } W_{\gamma}, \gamma \in \Gamma\}.
\]
We will take $e_{(\gamma,1)}$ as the highest weight vector $b^{\gamma}$. 
Define the set 
\[
 \Theta=\{(y,z) \in \bbr^2: 0 \leq y,z \leq 1, y^2+z^2=1\}.
\]
For $\gamma=(\gamma_1,\gamma_2) \in \Gamma$, define the  function 
\[
 g^{(\gamma_1,\gamma_2)} : \Theta \rightarrow \bbr
\]
sending $(y,z)$ to $y^{\gamma_1}z^{\gamma_2}$. Applying rotations on the co-ordinates appropriately, 
we get
\[
 \| b^{(\gamma_1,\gamma_2)}\|=\sup_{(y,z)\in \Theta } g^{(\gamma_1,\gamma_2)}(y,z). 
\]
Let us state one result  of \cite{Sau-2018aa} (see Proposition 3.2 in \cite{Sau-2018aa}).
\bppsn \label{cpt}
Let $\Theta$ be a compact subset of $\bbr^n$ and $f$ and $h$ are two real valued continuous   functions define on $\Theta$.
Let $x_0 \in \Theta$ be a point such that $|f(x_0)|=\|f\|= \sup_{x \in \Theta}|f(x)| \neq 0$ and $h(x_0) \neq 0$. Then one has
\[
 \frac{\|h^mf\|}{\|h^{m+1}f\|}\leq \frac{1}{|h(x_0)|}.
\]
\eppsn

\blmma \label{bound hwv}
Let $\epsilon_1=(1,0)$ and $\epsilon_2=(0,1)$. Then one has
\begin{enumerate} 
 \item
 $\sup_{\{\gamma \in \Gamma:\gamma_1=\gamma_2\}}\frac{\|b^{\gamma}\|}{\|b^{\gamma+\epsilon_1+\epsilon_2}\|} <\infty$.
  \item
 $\sup_{\{\gamma \in \Gamma:\gamma_1\geq \gamma_2\}}\frac{\|b^{\gamma}\|}{\|b^{\gamma+\epsilon_1}\|} <\infty$.
 \item
 $\sup_{\{\gamma \in \Gamma:\gamma_1\leq\gamma_2\}}\frac{\|b^{\gamma}\|}{\|b^{\gamma+\epsilon_2}\|} <\infty$. 
\end{enumerate}
\elmma
\prf It is not difficult to check that $g^{(1,1)}$ attains its maximum value at $y=1/\sqrt{2}$ and $z=1/\sqrt{2}$. Let $M=\sup_{(y,z)\in \Theta } yz$. 
\begin{enumerate}
 \item 
 If $\gamma_1=\gamma_2$ then $g^{(\gamma_1,\gamma_2)}=(yz)^{\gamma_1}$ and $g^{\gamma+\epsilon_1+\epsilon_2}=(yz)^{\gamma_1+1}$. Therefore,
 \[
  \sup_{\{\gamma \in \Gamma:\gamma_1=\gamma_2\}}\frac{\|b^{\gamma}\|}{\|b^{\gamma+\epsilon_1+\epsilon_2}\|} 
  = \sup_{\{\gamma \in \Gamma:\gamma_1=\gamma_2\}}\frac{\sup_{(y,z)\in \Theta}g^{(\gamma_1,\gamma_1)}}{\sup_{(y,z)\in \Theta}g^{(\gamma_1+1,\gamma_1+1)}}
  = \sup_{\{\gamma \in \Gamma:\gamma_1=\gamma_2\}}\frac{M^{\gamma_1}}{M^{\gamma_1+1}}=\frac{1}{M} <\infty.
 \]
 \item We have
 \begin{IEEEeqnarray}{rCl}
  \sup_{\{\gamma \in \Gamma:\gamma_1\geq \gamma_2\}}\frac{\|b^{\gamma}\|}{\|b^{\gamma+\epsilon_1}\|}&=& 
  \sup_{\{\gamma \in \Gamma:\gamma_1\geq \gamma_2\}}\frac{\sup_{(y,z)\in \Theta}g^{(\gamma_1,\gamma_2)}}{\sup_{(y,z)\in \Theta}g^{(\gamma_1+1,\gamma_2)}} \nonumber \\
  &=& \sup_{\{\gamma \in \Gamma:\gamma_1\geq \gamma_2\}}\frac{\sup_{(y,z)\in \Theta}y^{\gamma_1-\gamma_2}g^{(\gamma_2,\gamma_2)}}{\sup_{(y,z)\in \Theta}y^{\gamma_1-\gamma_2+1}g^{(\gamma_2,\gamma_2)}} \nonumber \\
  &=& \sqrt{2} < \infty \qquad \qquad \qquad \qquad \mbox{ (by Proposition } \ref{cpt}) \nonumber.
 \end{IEEEeqnarray}
 
\item We have
 \begin{IEEEeqnarray}{rCl}
  \sup_{\{\gamma \in \Gamma:\gamma_1\leq \gamma_2\}}\frac{\|b^{\gamma}\|}{\|b^{\gamma+\epsilon_2}\|}&=& 
  \sup_{\{\gamma \in \Gamma:\gamma_1\leq \gamma_2\}}\frac{\sup_{(y,z)\in \Theta}g^{(\gamma_1,\gamma_2)}}{\sup_{(y,z)\in \Theta}g^{(\gamma_1,\gamma_2+1)}} \nonumber \\
  &=& \sup_{\{\gamma \in \Gamma:\gamma_1\leq \gamma_2\}}\frac{\sup_{(y,z)\in \Theta}z^{\gamma_2-\gamma_1}g^{(\gamma_1,\gamma_1)}}{\sup_{(y,z)\in \Theta}z^{\gamma_2-\gamma_1+1}g^{(\gamma_2,\gamma_2)}} \nonumber \\
  &=& \sqrt{2} < \infty \qquad \qquad \qquad \qquad \mbox{ (by Proposition } \ref{cpt}) \nonumber.
 \end{IEEEeqnarray}

\end{enumerate}


\qed \\
Let $c>0$ be an upper bound in all the three inequalities  of the Lemma \ref{bound hwv}. Take $R=\{y,z,yz\}$. 
The following lemma says that $\mathcal{G}_R^c$ has a root $(0,0)$.
\blmma \label{path}
Let $\gamma \in \Gamma$. Then there is a path in  $\mathcal{G}_R^c$ joining $(0,0)$ and  $\gamma$ and of length  less than or equal to
$\max\{\gamma_1,\gamma_2\}$.
\elmma
\prf 
If $\gamma_1\geq\gamma_2$, then one possible path would be as follows.
\begin{IEEEeqnarray}{lCl}
(0,0)\rightarrow (1,1)\rightarrow (2,2)\rightarrow \cdots \rightarrow (\gamma_2,\gamma_2) \nonumber \\
\hspace{2in} \mbox{(by part(1) of the Lemma } \ref{bound hwv}) \nonumber \\
 (\gamma_2,\gamma_2) \rightarrow (\gamma_2+1,\gamma_2) \rightarrow \cdots \rightarrow
 (\gamma_1,\gamma_2)  \nonumber \\
  \hspace{2in} \mbox{ (by part(2) of the  Lemma } \ref{bound hwv}) \nonumber 
\end{IEEEeqnarray}
If $\gamma_1\leq\gamma_2$, then one possible path would be as follows.
\begin{IEEEeqnarray}{lCl}
(0,0)\rightarrow (1,1)\rightarrow (2,2)\rightarrow \cdots \rightarrow (\gamma_1,\gamma_1) \nonumber \\
\hspace{2in} \mbox{(by part(1) of the Lemma } \ref{bound hwv}) \nonumber \\
 (\gamma_1,\gamma_1) \rightarrow (\gamma_1+1,\gamma_1) \rightarrow \cdots \rightarrow
 (\gamma_1,\gamma_2)  \nonumber \\
  \hspace{2in} \mbox{ (by part(3) of the  Lemma } \ref{bound hwv}) \nonumber 
\end{IEEEeqnarray}
Moreover, the length of the path in first  case is $\gamma_1$ and in second case, it is $\gamma_2$. This settles the claim.
\qed 

\blmma \label{bounded leap}
For $1\leq m \leq n+1$, $l=n+1$ and $\gamma \in \Gamma$, one has 
\[
 u_m^le_{(\gamma,i)} \subset \mbox{span}\{ e_{(\beta,j)}: \beta \in \Gamma, \gamma_1 \leq \beta_1 \leq \gamma_1+1 \mbox { and }
  \gamma_{2} \leq \beta_2 \leq \gamma_{2}+1\}.
\] 
\elmma
\prf
Let $\epsilon_i=(0,0,\cdots, 0, \underbrace{1}_{i^{th}-\mbox{place}},0,\cdots 0)$. Then from equation ((13), page $210$, \cite{KliSch-1997aa}), we get 
\[
  u_{(1,0,\cdots ,0)}\otimes  u_{(\gamma_1,0,\cdots ,0,\gamma_{2})}= \oplus_{i=1}^n u_{(\gamma_1,0,\cdots ,0,\gamma_{2}) + e_i} 
\]                                                                                               
Hence $u_m^le_{(\gamma,i)}$ is in the span of matrix entries of the irreducible representations of  highest weight 
$(\delta_1, \delta_2, \cdots, \delta_{n+1})$ such that $\delta_1= \gamma_1$ or $\gamma_1+1$ and $\delta_{n+1}= \gamma_{2}$ or $\gamma_{2}+1$. 
Since  $u_m^le_{(\gamma,i)} \in \mathcal{O}(SU(n+1)/SU(n))$ and   
$\{e_{(\alpha,j)}: \alpha \in \Gamma, 1 \leq j \leq N_{\alpha}\}$ is a basis of $\mathcal{O}(SU(n+1)/SU(n))$,  we get the claim.
\qed 
\bthm \label{optimal}
Let $L_{(0,0)}$ be the Dirac operator  $e_{(\gamma,i)} \mapsto \max\{\gamma_1,\gamma_2\}e_{(\gamma,i)}$ 
acting on the Hilbert space $\clh_{\rho}$. Then  the triple
$(\mathcal{O}(SU(n+1)/SU(n)),\clh_{\rho},L_{(0,0)})$ is a $(2n+1)$-summable
equivariant spectral triple  of the  system $(C(SU(n+1)/SU(n)), SU(n+1),\tau)$.
The operator $L_{(0,0)}$ is optimal, i.e. if $D$ is any equivariant Dirac operator of
the $C^*$-dynamical system $(C(SU(n+1)/SU(n)), SU(n+1), \tau)$ acting on $\clh_{\rho}$ then there exist positive reals $a$ and $b$ such that
\[
 |D|\leq a|L_{(0,0)}|+b.
\]
\ethm 
\prf Clearly $L_{(0,0)}$ is a selfadjoint operator with compact resolvent.
That $L_{(0,0)}$ has bounded commutators with the generators $\left\{u_m^1 : m \in \{1,2,\cdots n+1\}\right\}$ of $\mathcal{O}(SU(n+1)/SU(n))$
follows from Lemma \ref{bounded leap}. This proves that the triple $(\mathcal{O}(SU(n+1)/SU(n)), \clh_{\rho},L_{(0,0)})$ is an 
equivariant spectral triple  of the  system $(C(SU(n+1)/SU(n)), SU(n+1), \tau)$.
From Weyl dimension formula, we have
\[
 N_{\gamma}=O(\gamma_1^{n-1}\gamma_2^{n-1}(\gamma_1+\gamma_2)).
\]
This shows that $L_{(0,0)}$ is $(2n+1)$-summable. Optimality follows from Proposition \ref{geq}.
\qed 
\bthm \label{spectral}
Spectral dimension of the  odd dimensional sphere $SU(n+1)/SU(n)$ is $2n+1$.
\ethm
\prf
It is a direct consequence of Theorem \ref{optimal}.
\qed
\newsection{$SO(2n+1)$ action on $\cls^{2n}$}
For $1 \leq i,j \leq 2n+1$, define a continuous map 
\begin{IEEEeqnarray}{rCl}
 u_j^i:SO(2n+1) & \rightarrow & \bbc;  \quad A \mapsto a_j^i \nonumber
\end{IEEEeqnarray} 
where $a_j^i$ is the $ij^{th}$ entry of $A \in SO(2n+1)$. 
The $C^*$-algebra  $C(SO(2n+1))$ is generated by elements of the set $\{u_{j}^{i}:1\leq i,j\leq 2n+1\}$.
In the same way,  define the generators  $\{v_{j}^{i}:1\leq i,j\leq 2n\}$  of $C(SO(2n))$. 
Define the map $\Phi: C(SO(2n+1)) \rightarrow C(SO(2n))$ 
as follows. 
\[ 
\Phi(u_{j}^{i})=\begin{cases} 
                             v_{j-1}^{i-1}, 
                             & \mbox{ if } i \neq 1 \mbox{ or } j \neq 1, \cr
			 \delta_{ij}, & \mbox{ otherwise. } \cr
			\end{cases}
\] 
The quotient space $SO(2n+1)/SO(2n)$ can be realized as the $2n$-dimensional sphere $\cls^{2n}$.  Also, each of 
the generators $\{u_{j}^{1}:   1\leq j\leq 2n+1\}$ can be viewed as 
projection on to a fixed real coordinate of a point in $\cls^{2n}\subset \bbr^{2n+1}$. The set of highest weights of all finite-dimensional irreducible co-representation of $C(SO(2n+1))$ 
can be described as follows.
\[
 \varLambda=\{(\lambda_1, \lambda_2, \cdots ,\lambda_{n}):\lambda_1\geq \lambda_2 \geq \cdots \geq \lambda_{n}\geq 0, \lambda_i's  \mbox { are  all either  integers or  half integers } \}.
\]
Using Zhelobenko branching rule (see Theorem 9, page 74, \cite{Zhe-1962aa} 
and Theorem 1.7 in \cite{Pod-1995aa}), we get  
\begin{displaymath}
I_{\lambda}= \begin{cases}
                            1 & \quad \mbox{if } \lambda_i=0 \quad \mbox{for all}\quad 2 \leq  i \leq n, \cr
                             0 & \quad \mbox{otherwise}. \cr
                            \end{cases}
\end{displaymath}  
The Hopf $*$-algebra $\mathcal{O}(SO(2n+1)/SO(2n))$ has a $\mathfrak{so}(2n+1)$ module structure induced by the pairing between  $\mathcal{O}(SO(2n+1)/SO(2n))$ and $\mathfrak{so}(2n+1)$ 
associated with the finite dimensional irreducible co-representation of highest weight $(1,0,\cdots,0)$ (similar to that defined in previous subsection). We 
call an element $b \in \mathcal{O}(SO(2n+1)/SO(2n))$ a highest weight vector  with highest weight $(\lambda_1,0, \cdots,0)$ if 
\[
 H_1(b)=2\lambda_1b, \quad  \quad H_i(b)= 0  \mbox{ for} \quad 2\leq i\leq n,
\]
and
\[
 E_i(b)=0 \qquad  \mbox{ for } 1\leq i \leq n.
\]
The following proposition describes a  highest weight vector for   highest weight $(\lambda_1,0, \cdots,0)$ explicitly.
\bppsn \label{new1}
Let $y= u_{2n+1}^{1}$. For $\gamma \in \bbn$, define $b^{\gamma}=y^{2\gamma}$. 
Then $b^{\gamma}$ is a  highest weight vector  in the algebra
$\mathcal{O}(SO(2n+1)/SO(2n))$ with highest weight $(\gamma,0, \cdots,0)$. 
\eppsn
\prf  It follows from a straightforward calculation.
\qed \\

Let $\Gamma=\{\gamma:\gamma \in \bbn\}$ and $R=\{y^2\}$. Define $e_{(\gamma,1)}=b^{\gamma}$. Since $\|b^{\gamma}\|=1$, one can show that for $c < 1$, one has $\gamma \rightarrow \gamma+1$ and hence the graph  
$\mathcal{G}_R^c$ has a root $0$. Further $\ell_0(\gamma)=\gamma$ and hence  associated length operator $L_0$ maps $e_{(\gamma,i)}$ to $\gamma e_{(\gamma,i)}$. To show that $L_0$ has bounded commutator 
with the generators of $C(SO(2n+1)/SO(2n))$, we need the following result.

\blmma \label{bounded leap1}
For $1\leq m \leq 2n+1$ and $l=1$, one has 
 \[
 u_m^le_{(\gamma,i)} \subset \mbox{span}\{ e_{(\beta,j)}: \beta \in \Gamma, \gamma-1 \leq \beta \leq \gamma+1 \}.
\] 
\elmma
\prf Proof follows by applying equation ((15), page $210$, \cite{KliSch-1997aa}) and taking similar steps as in Lemma \ref{bounded leap}.
\qed 
\bthm \label{optimal1}
Let $L_0$ be the Dirac operator  $e_{i}^{\gamma} \mapsto \gamma e_{i}^{\gamma}$ acting on the Hilbert space $\clh_{\rho}$. Then  the triple
$(\mathcal{O}(SO(2n+1)/SO(2n)),\clh_{\rho},L_0)$ is a $(2n)$-summable equivariant spectral triple  of the  system $(C(SO(2n+1)/SO(2n)), SO(2n+1),\tau)$.
The operator $L_0$ is optimal, i.e. if $D$ is any equivariant Dirac operator of the $C^*$-dynamical system $(C(SO(2n+1)/SO(2n)), SO(2n+1),\tau)$ acting on $\clh_{\rho}$ then there exist positive reals $a$ and $b$ such that
\[
 |D|\leq a|L_0|+b.
\]
\ethm 
\prf Clearly $L_0$ is a selfadjoint operator with compact resolvent.
That $L_0$ has bounded commutators with the generators $\left\{u_m^1 : m \in \{1,2,\cdots 2n+1\}\right\}$ of $\mathcal{O}(SO(2n+1)/SO(2n))$
follows from Lemma \ref{bounded leap1}. This proves that the triple $(\mathcal{O}(SO(2n+1)/SO(2n)),\clh_{\rho},L_0)$ is an 
equivariant spectral triple  of the  system $(C(SO(2n+1)/SO(2n)), SO(2n+1),\tau)$.
From Weyl dimension formula, we have
\[
 N_{\gamma}=O(\gamma_1^{2n-1}).
\]
This  shows that $L_0$ is $(2n)$-summable. Optimality follows from Proposition \ref{geq}.
\qed 
\bthm \label{spectral1}
Spectral dimension of the  even dimensional sphere $\cls^{2n}=SO(2n+1)/SO(2n)$ is $2n$.
\ethm
\prf It is a direct consequence of Theorem \ref{optimal1}.
\qed
\newsection{$SO(2n)$ action on $\cls^{2n-1}$}
For $1 \leq i,j \leq 2n$, define a continuous map 
\begin{IEEEeqnarray}{rCl}
 u_j^i:SO(2n) & \rightarrow & \bbc;  \quad A \mapsto a_j^i \nonumber
\end{IEEEeqnarray} 
where $a_j^i$ is the $ij^{th}$ entry of $A \in SO(2n)$. 
The $C^*$-algebra  $C(SO(2n))$ is generated by elements of the set
$\{u_{j}^{i}:1\leq i,j\leq 2n\}$. In the same way,  define the generators  $\{v_{j}^{i}:1\leq i,j\leq 2n-1\}$  of $C(SO(2n-1))$. 
Define the map $\Phi: C(SO(2n)) \rightarrow C(SO(2n-1))$ 
as follows. 
\[ 
\Phi(u_{j}^{i})=\begin{cases} 
                             v_{j-1}^{i-1}, 
                             & \mbox{ if } i \neq 1 \mbox{ or } j \neq 1, \cr
			 \delta_{ij}, & \mbox{ otherwise. } \cr
			\end{cases}
\] 
The quotient space $SO(2n)/SO(2n-1)$ can be realized as the $2n-1$-dimensional sphere $\cls^{2n-1}$.  Also, each of 
the generators $\{u_{j}^{1}:   1\leq j\leq 2n\}$ can be viewed as 
projection on to a fixed real coordinate of 
a point in $\cls^{2n-1}\subset \bbr^{2n}$. 
The set of highest weights of all finite-dimensional irreducible co-representation of $C(SO(2n))$ 
can be described as follows.
\[
 \varLambda=\{(\lambda_1, \lambda_2, \cdots ,\lambda_{n}):
   \lambda_1 \geq \cdots \geq\lambda_{n-1} \geq  |\lambda_{n}|, \lambda_n \in \bbn \cup \frac{\bbn}{2} \mbox{ and }
   \lambda_{i}-\lambda_{i+1} \in \bbn  \mbox{ for } 1 \leq i <n \}.
\]
Using Zhelobenko branching rule (see Theorem 10, page 77, \cite{Zhe-1962aa} 
and Theorem 1.7 in \cite{Pod-1995aa}), we get  
\begin{displaymath}
I_{\lambda}= \begin{cases}
                            1 & \quad \mbox{if } \lambda_i=0 \quad \mbox{for all}\quad 2 \leq  i \leq n, \cr
                             0 & \quad \mbox{otherwise}. \cr
                            \end{cases}
\end{displaymath}  
The Hopf $*$-algebra $\mathcal{O}(SO(2n)/SO(2n-1))$ has a 
$\mathfrak{so}(2n)$ module structure induced by the pairing between 
$\mathcal{O}(SO(2n)/SO(2n-1))$ and $\mathfrak{so}(2n)$ 
associated with the finite dimensional irreducible co-representation of 
highest weight $(1,0,\cdots,0)$ (similar to that defined in previous subsection). We 
call an element $b \in \mathcal{O}(SO(2n)/SO(2n-1))$ a highest weight vector  with highest weight $(\lambda_1,0, \cdots,0)$ if 
\[
 H_1(b)=\lambda_1b, \quad  \quad H_i(b)= 0  \mbox{ for} \quad 2\leq i\leq n,
\]
and
\[
 E_i(b)=0 \qquad  \mbox{ for } 1\leq i \leq n.
\]
The following proposition describes a  highest weight vector for   highest weight $(\lambda_1,0, \cdots,0)$ explicitly.
\bppsn \label{new2}
Let $y= u_{2n+1}^{1}$. For $\gamma \in \bbn$, define $b^{\gamma}=y^{\gamma}$. 
Then $b^{\gamma}$ is a  highest weight vector  in the algebra
$\mathcal{O}(SO(2n)/SO(2n-1))$ with highest weight $(\gamma,0, \cdots,0)$. 
\eppsn
\prf  It follows from a straightforward calculation.
\qed \\

Let $\Gamma=\{\gamma:\gamma \in \bbn\}$ and $R=\{y\}$. Define $e_{(\gamma,1)}=b^{\gamma}$. Since $\|b^{\gamma}\|=1$, one can show that for $c < 1$, one has $\gamma \rightarrow \gamma+1$ and hence the graph  
$\mathcal{G}_R^c$ has a root $0$. Moreover $\ell_0(\gamma)=\gamma$ and hence  associated length operator $L_0$ maps $e_{(\gamma,i)}$ to $\gamma e_{(\gamma,i)}$. To show that $L_0$ has bounded commutator 
with the generators of $C(SO(2n)/SO(2n-1))$, we need the following result.

\blmma \label{bounded leap2}
For $1\leq m \leq 2n$ and $l=1$, one has 
\[
 u_m^lu_{i}^{\gamma} \subset \mbox{span}\{ u_{i}^{\beta}: \gamma_1-1 \leq \beta_1 \leq \gamma_1+1 \}
\]
\elmma
\prf Proof follows by applying equation ((14), page $210$, \cite{KliSch-1997aa}) and taking similar steps as in Lemma \ref{bounded leap}.
\qed 
\bthm \label{optimal2}
Let $L_0$ be the Dirac operator  $e_{i}^{\gamma} \mapsto \gamma e_{i}^{\gamma}$ acting on the Hilbert space $\clh_{\rho}$. Then  the triple
$(\mathcal{O}(SO(2n)/SO(2n-1)),\clh_{\rho},L_0)$ is a $(2n-1)$-summable equivariant spectral triple  of the  system $(C(SO(2n)/SO(2n-1)), SO(2n),\tau)$.
The operator $L_0$ is optimal, i.e. if $D$ is any equivariant Dirac operator of the $C^*$-dynamical system 
$(C(SO(2n)/SO(2n-1)), SO(2n),\tau)$ acting on $\clh_{\rho}$ then there exist positive reals $a$ and $b$ such that
\[
 |D|\leq a|L_0|+b.
\]
\ethm 
\prf Clearly $L_0$ is a selfadjoint operator with compact resolvent.
That $L_0$ has bounded commutators with the generators $\left\{u_m^1 : m \in \{1,2,\cdots 2n\}\right\}$ of $\mathcal{O}(SO(2n)/SO(2n-1))$
follows from Lemma \ref{bounded leap2}. This proves that the triple $(\mathcal{O}(SO(2n)/SO(2n-1)),\clh_{\rho},L_0)$ is an 
equivariant spectral triple  of the  system $(C(SO(2n)/SO(2n-1)), SO(2n),\tau)$.
From Weyl dimension formula, we have
\[
 N_{\gamma}=O(\gamma_1^{2n-2}).
\]
This  shows that $L_0$ is $(2n-1)$-summable. Optimality follows from Proposition \ref{geq}.
\qed 
\bthm \label{spectral2}
Spectral dimension of the  odd dimensional sphere $SO(2n)/SO(2n-1)$ is $2n-1$.
\ethm
\prf It is a direct consequence of Theorem \ref{optimal2}.
\qed

\noindent{\sc Bipul Saurabh} (\texttt{saurabhbipul2@gmail.com})\\
         {\footnotesize Harish Chandra Research Institute , Chhatnag Road, Jhunsi,
Allahabad 211019,  INDIA}

\end{document}